\documentclass[11pt]{amsart}
\usepackage{amsmath}
\usepackage{epsfig}
\usepackage{graphicx}
\usepackage{eufrak}
\usepackage{dsfont}
\usepackage[english]{babel}
\usepackage{color}
\usepackage{tikz-cd}

\usepackage{palatino, mathpazo}
\usepackage{amscd}      
\usepackage{amssymb}
\usepackage{dsfont}
\usepackage{xypic}      
\LaTeXdiagrams          
\usepackage[all,v2]{xy}
\xyoption{2cell} \UseAllTwocells \xyoption{frame} \CompileMatrices
\usepackage{latexsym}
\usepackage{epsfig}
\usepackage{enumerate}

\usepackage{latexsym}
\usepackage{epsfig}
\usepackage{amsfonts}
\usepackage{enumerate}
\usepackage{mathrsfs}

\newbox\mybox
\def\overtag#1#2#3{\setbox\mybox\hbox{$#1$}\hbox to
  0pt{\vbox to 0pt{\vglue-#3\vglue-\ht\mybox\hbox to \wd\mybox
      {\hss$\ss#2$\hss}\vss}\hss}\box\mybox}
\def\undertag#1#2#3{\setbox\mybox\hbox{$#1$}\hbox to 0pt{\vbox to
    0pt{\vglue#3\vglue\ht\mybox\hbox to \wd\mybox
      {\hss$\ss#2$\hss}\vss}\hss}\box\mybox}
\def\lefttag#1#2#3{\hbox to 0pt{\vbox to 0pt{\vss\hbox to
      0pt{\hss$\ss#2$\hskip#3}\vss}}#1}
\def\righttag#1#2#3{\hbox to 0pt{\vbox to 0pt{\vss\hbox to
      0pt{\hskip#3$\ss#2$\hss}\vss}}#1}
\let\ss\scriptstyle

\def\Dot{\lower.2pc\hbox to 2pt{\hss$\bullet$\hss}}
\def\Circ{\lower.2pc\hbox to 2pt{\hss$\circ$\hss}}
\def\Vdots{\raise5pt\hbox{$\vdots$}}
\newcommand\lineto{\ar@{-}}
\newcommand\dashto{\ar@{--}}
\newcommand\dotto{\ar@{.}}

\allowdisplaybreaks[3]

\newtheorem{prop}{Proposition}[section]
\newtheorem{lem}[prop]{Lemma}

\newtheorem{cor}[prop]{Corollary}
\newtheorem{thm}[prop]{Theorem}
\newtheorem{rmk}[prop]{Remark}
\newtheorem{example}{Example}

\newtheorem{defn}[prop]{Definition}

            {\nolinebreak $\Box$ \end{trivlist}}

\newcommand{\noprint}[1]{}

\renewcommand{\tilde}{\widetilde}

\newcommand{\Ext}{\mbox{Ext}}

\newcommand{\Hom}{\mbox{Hom}}

\newcommand{\scM}{\mathscr{M}}

\newcommand{\XX}{{\mathfrak X}}

\renewcommand{\SS}{{\mathfrak S}}

\newcommand{\ZZ}{{\mathfrak Z}}

\newcommand{\Mm}{{\mathfrak m}}

\newcommand{\zz}{{\mathbb Z}}

\newcommand{\aaa}{{\mathbb A}}
\newcommand{\ttt}{{\mathbb T}}

\renewcommand{\ll}{{\mathbb L}}
\newcommand{\qq}{{\mathbb Q}}
\newcommand{\pp}{{\mathbb P}}

\newcommand{\sE}{{\mathcal E}}

\newcommand{\sO}{{\mathcal O}}
\newcommand{\sR}{{\mathcal R}}

\newcommand{\rM}{\mathscr{M}}

\newcommand{\Oo}{{\mathscr O}}
\newcommand{\Rr}{{\mathscr R}}

\DeclareMathOperator{\lci}{lci}

\DeclareMathOperator{\coker}{coker}

\DeclareMathOperator{\mult}{mult}
\DeclareMathOperator{\El}{El}
\DeclareMathOperator{\No}{No}
\DeclareMathOperator{\Cu}{Cu}
\DeclareMathOperator{\Ta}{Ta}
\DeclareMathOperator{\Tr}{Tr}
\DeclareMathOperator{\GL}{GL}

\newcommand{\ob}{\mathop{\rm ob}}

\newcommand{\spec}{\mathop{\rm Spec}\nolimits}

\numberwithin{equation}{subsection}

\newcommand {\mat}      [1] {\left(\begin{array}{#1}}
\newcommand {\rix}          {\end{array}\right)}

\title[Deformation of minimally elliptic  singularities]{Equivariant deformation of minimally elliptic  singularities}
\date{January, 2025}

\author{Sagnik Das}
\address{Department of Mathematics\\ University of Kansas\\ 405 Snow Hall 1460 Jayhawk Blvd\\Lawrence KS 66045 USA} 
\email{das.sagnik@ku.edu}

\author{Yunfeng Jiang}
\address{Department of Mathematics\\ University of Kansas\\ 405 Snow Hall 1460 Jayhawk Blvd\\Lawrence KS 66045 USA} 
\email{y.jiang@ku.edu}

\begin{document}
\sloppy \maketitle
\begin{abstract}
We study certain equivariant deformation components of minimally elliptic surface singularities under finite group actions. Interesting examples include cyclic quotients of simple elliptic singularities and finite group quotients of cusp singularities, where the resulting quotients remain simple elliptic and cusp singularities, respectively. In cases where the minimally elliptic singularities are locally complete intersection (lci) singularities, we identify equivariant deformation components of general type surfaces containing such singularities that admit a perfect obstruction theory. 
\end{abstract}

\maketitle

\tableofcontents

\section{Introduction}

Minimally elliptic singularities, introduced by Laufer in \cite{Laufer}, are a very important class of normal surface singularities in both singularity theory and birational geometry. Simple elliptic and cusp singularities are examples of minimally elliptic singularities; and they also appear as log canonical (lc) singularities in the KSBA compactification of surfaces of log general type. Their smoothings have been extensively studied;  for example, \cite{Looijenga}, \cite{Looijenga2}, \cite{Wahl_Ann}, \cite{Wahl_Duke}, \cite{Wahl_Math-Ann}, \cite{GHK15}, and  \cite{Engel}.

Motivated by the construction of the virtual fundamental class on the KSBA moduli space of surfaces of general type in \cite{Jiang_2022}, we study the equivariant smoothing and deformation of locally complete intersection (lci) minimally elliptic singularities under a finite group action. The equivariant smoothing of simple elliptic and cusp singularities was previously studied in \cite{Jiang_2023}, \cite{Jiang_cusp}. The main goal of this paper is to prove $G$-equivariant versions of the deformation results for minimally elliptic singularities found in  \cite{Wahl_Ann}, \cite{Wahl_Duke}, \cite{Wahl_Math-Ann}.

We prove that every equivariant deformation component of an slc surface, originating from the $G$-equivariant deformation of an lci minimally elliptic singularity, admits a perfect obstruction theory. Furthermore, we generalize a result of J. Wahl on $G$-equivariant deformation components that can deform to lower-degree minimally elliptic singularities via partial resolutions.
 Throughout this paper, we work over an algebraically closed field $\mathbf{k}$.

\subsection{Background}

Let $(S,0)=(\spec(R),0)$ be a normal surface singularity.  Let $G$ be a finite group acting effectively on the surface $(S,0)$.  This means that $G$ acts on the integral domain $R$ and the quotient germ singularity $\overline{S}=(S/G,0)=(\spec(R^G), 0)$ is also a normal singularity, since $R^G$ is also a normal integral domain. 

For the normal surface singularity $(\spec(R),0)$, the deformation space $D_R$ is defined as 
$$
D_R=\{\text{deformations~} \spec(\sR)\to \spec(A)\}/\text{equivalence}
$$
where $A$ is a local Artin $\mathbf{k}$-algebra.  We let 
$D_R^G\subset D_R$ be the subspace 
consisting of $G$-equivariant deformations $\spec(\sR)\to \spec(A)$.
Let $X\to \spec(R)$ be the minimal resolution with exceptional divisors $E_1, \cdots, E_r$.  Then the action of $G$ lifts to an action on $X$ such that $X\to \spec(R)$ is $G$-equivariant.  The quotient $X/G\to \spec(R^G)$ is a partial resolution of $\spec(R^G)$. 

Let $\XX\to \spec(A)$ be the $G$-equivariant deformation of  $X\to \spec(R)$  such that the $G$-equivariant deformation space $D_X^G$ is subspace of $D_X$, and 
the quotient $\XX/G\to \spec(A)$ is a deformation of $X/G\to \spec(R^G)$.
We first have the following result which generalizes one of the main results in \cite{Wahl_Ann}.

\begin{thm}\label{thm_Wahl-1}(Theorem \ref{cor_J.Wahl1})
   Given $h^2(\sO_{[X/G]})= 0$,  the  $G$-equivariant deformation  $\XX\to \spec(A)$ blows down to give a $G$-equivariant deformation of $\spec(R)$ if and only if $h^1(\sO_{[X/G]})$ remains constant. 
\end{thm}

We let $B^G\subset D_X^G$ be the subspace consisting of the  $G$-equivariant deformations  $\XX\to \spec(A)$ that blow down to give $G$-equivariant deformations of $\spec(R)$ in Theorem \ref{thm_Wahl-1}. 

We have several remarks. Let $B\subset D_X$ be the without considering the $G$ action. In the cases that $(\spec(R),0)$ are A, D, E type (RDP) singularities, then there is a  morphism $\Phi: B\to D_R$
such that $\Phi(B)$ is an irreducible component $A$ called the Artin component. 
The morphism $\Phi: B\to A$ is Galois with Weyl group $W$ coming from the graph of $R$.
In the case that $(\spec(R),0)$ is a rational singularity, $B=D_X$ is smooth since $h^1(\sO_X)=0$.

We are interested in the case that $(S,0)$ is Gorenstein.   The singularity $(S,0)$ is Gorenstein if and only if the dualizing module $\omega_R$ is invertible.  The quotient $\overline{S}$ is $\qq$-Gorenstein since $G$ is finite. 
We are interested in the case that the quotient $\overline{S}$ is still Gorenstein.
The quotient $\overline{S}$ is  Gorenstein if and only the dualizing sheaf 
$\omega_{S}^G=\omega_{\overline{S}}$ is still 
invertible. 

\subsection{Minimally elliptic singularity}

Consider the condition $h^1(\sO_X)=1$, in which case the singularity is called a minimally elliptic singularity from \cite{Laufer}.
Simple elliptic singularities and cusp singularities are examples of minimally elliptic singularities, which are the index one lc singularities in the KSBA compactification of log general type surfaces.  There exists cyclic group actions on such singularities which can have index one or higher.  If the quotient has index bigger than one, then it is $\qq$-Gorenstein. For instance, the $\zz_2$-quotient of a cusp singularity $(S,0)$ that has a resolution graph of a tree of rational curves is a rational singularity, and $\qq$-Gorenstein.

For the $G$-equivariant minimal resolution 
$X\to (\spec(R),0)$, let $Z$ be the fundamental cycle. 
Let $R_Z^G\subset D_X^G$ be the subspace consisting of $G$-equivariant deformations such that the exceptional cycle $Z$ lifts $G$-equivariantly.

\begin{thm}\label{thm_Wahl-2}
    In the case where $(\spec(R),0)$ is a minimally elliptic singularity together with a $G$-action, we have $$ R_Z^G = B^G $$ and the map $B^G\to D_R^G$ factors through $B^G\to B^G/W\to NF^G\to D_R^G$, where $W$ is the Weyl group, and $NF^G$ is the subspace in $D_R^G$ of normal flat deformations equivariant under the $G$-action.

    If $\mult(R)\ge 3$, $B^G/W\cong NF^G$. 

    If $\mult(R)\ge 10$, $NF^G$ consists of equivariant irreducible components of $D_R^G$.

    In all the above cases, $B^G$ give rise to the deformation subspace of $(\spec(R^G),0)$ induced from the $G$-equivariant deformations of 
    the minimal resolution of $(\spec(R),0)$.
\end{thm}

Suppose $R^G$ is still Gorenstein, and we have 
\[
\xymatrix{
X\ar[r]\ar[d]&\spec(R)\ar[d]\\
\overline{X}\ar[r]&\spec(R^G)
}
\]
and $h^1(\sO_X)$ is $G$-invariant. 

\begin{cor}\label{cor_deform_compo_G_intro}
    The space $B^G$ gives the deformation space in $D_{\overline{X}}\subset D_{R^G}$ induced from the $G$-equivariant deformation of $R$ through $X\to \spec(R)$ blowing down. 
\end{cor}

The above results can apply to the $G$-equivariant deformation of special type of minimally elliptic singularities--simple elliptic, cusp, tacnode, triangle singularities.  
Let $\El, \Cu, \No, \Ta, \Tr$ represent the simple elliptic, cusp, Node, Tacnode, and Triangle singularities. 
We have 

\begin{thm}\label{thm_G_equisingular_minimal_intro}(Theorem \ref{thm_G_equisingular_minimal})
    Let $(\spec(R),0)$ be a minimally elliptic singularity and $X\to \spec(R)$ the $G$-equivariant minimal resolution. Let $Z$ be the fundamental cycle and let $d:=-Z\cdot Z$ which is given by the following formula for the type of minimally elliptic singularities.
    \begin{equation}\label{eqn_ZZ_G_Cu_El}
    \begin{array}{ll}
      \Cu(d_1, d_2,\cdots, d_r):   & -Z\cdot Z=\begin{cases}
          \sum (d_i-2) & r>1\\
          d & r=1.
      \end{cases}  \\
      \No(d): & -Z\cdot Z=d \\
    \Ta(d_1,d_2): & -Z\cdot Z=d_1+d_2-4 \\
    \Tr(d_1,d_2,d_3): & -Z\cdot Z=d_1+d_2+d_3-6
    \end{array}
\end{equation}
    Then the $G$-equivariant deformations  (which are $-Z\cdot Z$-constant deformations) form an irreducible component of the deformation space $D_R^G$. 
\end{thm}

\begin{thm}\label{thm_G_singular_minimal_deformation_intro}(Theorem \ref{thm_G_singular_minimal_deformation})
    Let $(\spec(R),0)$ be a minimally elliptic singularity and $X\to \spec(R)$ the $G$-equivariant minimal resolution. Let $Z$ be the fundamental cycle and let $d:=-Z\cdot Z$.

    Let $\Upsilon$ be a configuration of a rational singularity given by 
\begin{equation}\label{eqn_configuration_rational}
\xymatrix@R=7pt@C=24pt@M=0pt@W=0pt@H=0pt{
  &\undertag{\bullet}{-b_1}{1pt}\lineto[r]
  &\undertag{\bullet}{-b_2}{1pt}\dashto[rr]&&
  \undertag{\bullet}{-b_t}{1pt}&&}
  \end{equation} \\ 
in Proposition \ref{prop_Wahl_singulairty_G}.
Then we have the following $G$-equivariant deformations 
    \begin{equation}\label{eqn_deformations_G_Cu_El-2}
    \begin{array}{ll}
      (1)   &  \Cu(d_1,b_1,\cdots, b_t, d_2,\cdots, d_r)\longrightarrow \Cu(d_1+d_2-1,d_3,\cdots, d_r) (r>2) \\
       (2) & \Cu(d_1,b_1,\cdots, b_t, d_2)\longrightarrow \Cu(d_1+d_2-3)\\
       (3)  & \Cu(d, b_1,\cdots, b_t)\longrightarrow \El(d-1).
    \end{array}
\end{equation}
\end{thm}

\subsection{Application to perfect obstruction theory}\label{subsec_POT_intro}

We consider a special case of the deformation space $B^G$, which we call the moduli space.  Suppose that the singularity $(\spec(R),0)$ is lci such that the quotient $(\spec(R^G),0)$ is a minimally elliptic singularity, then from Corollary \ref{cor_deform_compo_G_intro}, the moduli space $B^G$ gives the $G$-equivariant deformation component of $\spec(R^G)$.

We are interested in the deformation components of semi-log-canonical (slc) surfaces containing minimally elliptic singularities.
Suppose that we have slc surfaces of \cite{Kollar-Shepherd-Barron} such that they contain only such singularities $\spec(R^G)$.  We further assume that there is a KSBA moduli space 
$M_B^G$ of deformation families for such  slc surfaces.  For each such a singularity, we take the local lci covering Deligne-Mumford stack which is given by 
$[\spec(R)/G]$, and let $M_B^{G,\lci}\to M_B^G$ be the moduli stack of lci covers constructed in \cite{Jiang_2022}. This moduli stack $M_B^{G,\lci}$ is given by  the $G$-equivariant deformation component of the lci cover $\SS^{\lci}\to S$, where $S$ is an slc surface containing the singularity $\spec(R^G)$, and $\SS^{\lci}$ is the lci cover Deligne-Mumford stack.  The deformation component $B^G\subseteq M_B^{G,\lci}$.
Let 
$$
p:  \SS^{\lci}\to M_B^{G,\lci}
$$
be the universal family of the moduli stack.  

\begin{cor}\label{cor_POT_intro}
The moduli stack of lci covers  $M_B^{G,\lci}$ admits a perfect obstruction theory.  
\end{cor}
\begin{proof}
     For the universal family 
       $p: \SS^{\lci}\to M_B^{G,\lci}$, let  $\ll_{\SS^{\lci}/M_B^{G,\lci}}$ be the relative cotangent complex and $\ll_{M_B^{G,\lci}}$ be the cotangent complex of $M_B^{G,\lci}$. Then the Kodaira-Spencer map defines a morphism 
       $$
       \phi: Rp_*(\ll_{\SS^{\lci}/M_B^{G,\lci}}\otimes \omega_{\SS^{\lci}/M_B^{G,\lci}}^{\bullet})[-1]\to \ll_{M_B^{G,\lci}}
       $$
       where $\omega_{\SS^{\lci}/M_B^{G,\lci}}^{\bullet}$ is the relative dualizing complex of the universal family which is a line bundle due to the fact that the  family is lci. Thus, the main result in \cite{Jiang_2022} implies that $\phi$ is a perfect obstruction theory.
\end{proof}

Corollary \ref{cor_POT_intro} provides strong evidence for the perfect obstruction theory on the moduli stack of lci covers constructed in \cite{Jiang_2022}.

\subsection*{Outline}

The paper is outlined as follows. In Section \ref{sec_G_equisingular_deformation}, we prove that the $G$-equivariant deformation of the minimal resolution of a surface singularity blows down if and only if the first cohomology of the structure sheaf remains constant under the $G$-action. 
In 
Section \ref{sec_simultaneous_G}, we study the above result 
for minimally elliptic singularities. 
In Section 
\ref{sec_deformation_elliptic_singularity}, we generalize the result of Wahl for the equivariant deformation of minimally elliptic singularities which involves ranging the degrees of the fundamental cycle. 

\subsection*{Acknowledgments}

Y. J. thanks Professor J. Wahl for his email correspondence of surface singularities. 
This work is partially supported by NSF DMS-2401484 and Simons Collaboration Grant.

\section{$G$-equivariant Equisingular deformation of surface singularities}\label{sec_G_equisingular_deformation}
In this section we generalize Section 0 and some portion of Section 1 of \cite{Wahl_Ann} to  quotient stacks $[X/G]$, where $G$ is a finite group acting effectively on a scheme $X$. First we recall some definitions and observations for this case. \\
\begin{defn}
    (\cite{Jacob_equi}) Let $X$ be a smooth scheme over field $\mathbf{k}$ and $G$ is a finite group acting on $X$ effectively. Let $\mathscr{C}$ be the category of Artin $\mathbf{k}$-algebras and $A\in \mathscr{C}$. A $G$-equivariant deformation of $X$ to $A$ is a commutative diagram 
    \begin{equation}\label{eqn:Gamma_bdle}
\xymatrix{
X\ar[d]\ar[r] & \XX \ar[d]\\
\spec(\mathbf{k})\ar[r] & \spec(A),
}    
\end{equation}
where $\XX$ is flat and separated over $\spec(A)$ where $G$ acts on $\XX$ effectively and the induced map $X \longrightarrow \XX {\times}_A \mathbf{k}$ is a $G$-equivariant isomorphism. 

Two such deformations $\XX$, $\XX$' are isomorphic if there exist a $G$-equivariant isomorphism $\phi:\XX \rightarrow \XX'$ over $\spec(A)$.
\end{defn}
Now we consider the following situation where we have a $G$-equivariant deformation as above, which gives a diagram 
\begin{equation}\label{eqn:Gamma_bdle}
\xymatrix{
[X/G]\ar[d]\ar[r] & [\XX/G] \ar[d]\\
\spec(\mathbf{k})\ar[r] & \spec(A).
}    
\end{equation}
The stack  $[X/G]$ is isomorphic to $[\XX/G] {\times}_A \mathbf{k}$, hence we have a diagram for deformation of $[X/G]$. Now $\XX$ is flat over $\spec(A)$,  which implies  $[\XX/G]$ is flat over $\spec(A)$. The action of $G$ on $\XX$ is proper which implies $[\XX/G]$ is separated over $\spec(A)$. We define the \v{C}ech cohomology groups for $[\XX/G]$. 
\begin{defn}
    \cite{Alper_moduli} Take an \'{e}tale covering $\mathscr{U}= \{U_i \rightarrow [\XX/G]\}_{i \in I}$. Now the \v{C}ech complex of $[\XX/G]$ with respect to the covering $\mathscr{U}$ is the complex $\check{\mathrm{C}}^{\bullet} (\mathscr{U}, {\sO}_{[\XX/G]})$ with the cohomology groups \\
     $\check{\mathrm{C}}^{n} (\mathscr{U}, {\sO}_{[\XX/G]}) = {\prod}_{(i_0,i_1,..,i_n)\in I^{n+1}}\sO_{[\XX/G]}(U_{i_0} {\times}_{[\XX/G]} U_{i_1} {\times}_{[\XX/G]} .. {\times}_{[\XX/G]} U_{i_n}) $ and differentials \\
     $$d^n : \check{\mathrm{C}}^{n} (\mathscr{U}, {\sO}_{[\XX/G]}) \rightarrow \check{\mathrm{C}}^{n+1} (\mathscr{U}, {\sO}_{[\XX/G]})$$ 
     defined by 
     $$(s_{i_0 i_1 .. i_n}) \longmapsto ({\sum}^{n+1}_{k=0} (-1)^k {p}^{*}_{\hat{k}} (s_{i_0..\hat{i}_k .. i_n}))_{(i_0, i_1,.., i_n)}  $$ 
     where ${p}^{*}_{\hat{k}}$ is the standard "Forgetting the $k$-th component map".  Now the \v{C}ech cohomology groups of $\sO_{[\XX/G]}$ with respect to the cover $\mathscr{U}$ is defined as the cohomology of the above complex and is denoted by $\check{\mathrm{H}}^n(\mathscr{U}, {\sO}_{[\XX/G]})$ for all $n$.
\end{defn}

 \begin{rmk}
         In place of the structure sheaf the same construction gives cohomology groups for any abelian sheaf on $[\XX/G]$. 
     \end{rmk}

Now $\XX$ is a separated $A$- scheme which implies it is a scheme with affine diagonal. Hence by \cite{Alper_moduli}, $[\XX/G]$ is also a stack with affine diagonal. This stack is also DM stack because $G$ is smooth affine finite group scheme over $\mathbb{C}$. 
We have the following proposition for DM stack with affine diagonal. 
\begin{prop}
    (\cite[Lemma 6.1.27, Proposition 6.1.28]{Alper_moduli}) Let $\XX$ be a DM stack with affine diagonal. If $\mathscr{U}$ is an affine \'{e}tale covering of $\XX$, then $H^i(\XX, F) = \check{\mathrm{H}}^n(\mathscr{U},F)$ for all quasi-coherent sheaf on $\XX$. 
\end{prop}
\begin{rmk}
    In the above proposition, the cohomology in the left is the cohomology groups defined by the i-th right derived functor of the global section functor $\Gamma : QCoh(\XX) \rightarrow Ab$. 
\end{rmk}
 The stack $[\XX/G]$ is flat over $A$, hence the \v{C}ech complex of the structure sheaf of it is a complex of flat modules over $A$. So by Artinian principle of exchange described in \cite{Wahl_Ann} we have the following proposition. 
 \begin{prop}\label{prop_0.4_Wahl}
     If $X, \XX, G$ are all defined as above, then we have natural maps \\
     $\Phi^i_A : H^{i} ([\XX/G]) {\bigotimes}_A k \rightarrow H^{i}([X/G])$ for all $i\geq0$, which satisfies the following properties. \\
     (a) If $\Phi^i_A$ is surjective then it is an isomorphism. \\
     (b) If $H^{i} ([\XX/G])$ is $A$-flat then $\Phi^i_A$  is injective.\\ 
     (c) If $\Phi^i_A$ is injective and $\Phi^{i-1}_A$ is an isomorphism then $\Phi^i_A$  is an isomorphism. \\
     (d) Any two of the following  will imply the third
     \begin{enumerate}
         \item $\Phi^i_A$ is an isomorphism. 
         \item $\Phi^{i-1}_A$ is an isomorphism. 
         \item $H^{i} ([\XX/G])$ is $A$- flat. 
     \end{enumerate}
     \end{prop}
     Now we talk about the equivariant blow-down deformations. Let $X, \spec(R)$ (dim$R \geq 2$), be two integral normal schemes and $G$ is a finite group as above acting effectively on both. Now let $f: X \rightarrow 
     \spec(R)$ be a $G$-equivariant proper, birational map. This induces a proper, birartional map of quotient DM stacks $\overline{f}: [X/G] \rightarrow [\spec(R)/G]$ \cite{Kresch_Tsch}. Here we observe dim $[\spec(R)/G]$ = dim$R \geq 2$ and the course moduli spaces (which are also the quotients $X/G$ and $\spec(R^G)$) are integral, normal and dim $(R^G) \geq 2$. So we have $R^G = \Gamma(X/G) = \Gamma([X/G])$.  
     Generalizing  \cite[Lemma 1.2]{Wahl_Ann} we have
     \begin{lem}
         For $A \in \mathscr{C}$, let us have a commutative diagram of $G$-equivariant deformations of $X, \spec(R)$ and $f$ over $\spec(A)$ as below

       \begin{equation}
           \begin{tikzcd}
X \arrow[rd] \arrow[dd, dotted] \arrow[r] & \XX \arrow[rd] \arrow[dd, dotted] &                               \\
                                          & \spec(R) \arrow[ld] \arrow[r]      & \spec(\overline{R}) \arrow[ld] \\
\spec(\mathbf{k}) \arrow[r]                         & \spec(A)                           &                              
\end{tikzcd}
       \end{equation}  
 Then $\overline{R}^G = \Gamma ([\XX/G])$.         
     \end{lem}
     \begin{proof}
         Descending the diagram to quotients, it is again a special case of lemma 1.2 of \cite{Wahl_Ann}. We observe that the $G$-equivariant blow down of $\XX$ is equivalent to  flatness of $\Gamma([\XX/G])$ or surjectivity of $\Gamma([\XX/G]) \rightarrow \Gamma([X/G])$. 
     \end{proof}
     Looking at the observation in the above proof we have the following result. 

     \begin{thm}\label{cor_J.Wahl1}
         Let $X, \spec(R), f$ be all defined as above and   $\XX$ over $\spec(A)$ is a $G$-equivariant deformation of $X$. Then
         \begin{enumerate}
         \item If $H^s(\sO_{[X/G]})=0$ for some $s\geq 3$ and if  $H^i(\sO_{[\XX/G]})$ is $A$-flat for $1 \leq i \leq s-1$, then $\XX$ blows down equivariantly. \\
         \item  If $H^2(\sO_{[X/G]})=0$, then $\XX$ blows down equivariantly iff $H^1(\sO_{[\XX/G]})$ is $A$-flat.\\
         \item  If $H^1(\sO_{[X/G]})=0$, then it always blows down equivariantly. 
         \end{enumerate}
     \end{thm}
     \begin{proof}
         This is the equivariant version of \cite[Theorem 1.4]{Wahl_Ann}, which is from Proposition \ref{prop_0.4_Wahl}.
     \end{proof}

\section{Simultaneous resolution of normal surface singularities under $G$-action}\label{sec_simultaneous_G}

\subsection{Some preliminaries}\label{subsec_preliminary}

Let us review the lifting of divisors. Let $X\to \spec(R)$ be the minimal resolution with exceptional divisors $E_1, \cdots, E_r$.  Then the action of $G$ lifts to an action on $X$ such that $X\to \spec(R)$ is $G$-equivariant.  Suppose we have a $G$-equivariant deformation 
$$
\XX\to \spec(\Oo)
$$
over a complete discrete ring such that the generic fiber is 
$\XX_K\to \spec(K)$ with $K$ the quotient field of $\Oo$. 
Let $\ZZ_K\subset \XX_K$ be a proper divisor over $K$.  Since $\XX$ is 
non-singular,  $\ZZ_K$ lifts to 
$\ZZ\subset \XX$, flat and proper over $\Oo$, and specializes to an effective exceptional divisor $Z$ on $X$.  We say $Z$ lifts to $\XX$ in this case. 

The standard semi-continuity theorem gives 
$$\chi(\sO_{\ZZ_K}(n\ZZ_K))=\chi(\sO_Z(nZ)); \quad  \ZZ_K\cdot \ZZ_K=Z\cdot Z.$$
If $\ZZ_K$ is reduced and irreducible, $\chi(\sO_{\ZZ_K})=1-h^1(\sO_{\ZZ_K})\le 1$ which implies that $\chi(\sO_Z\le 1)$.  Actually \cite[Lemma 1.2]{Wahl_Duke} proved that such exceptional divisors $Z$ must be finite. 

The following proposition is obviously true $G$-equivariantly.

\begin{prop}\label{prop_blowing_down}(\cite[Proposition 1.6]{Wahl_Duke})
 If the deformation $\XX\to \spec(\Oo)$ blows down to give a deforamtion 
 $\spec(\Rr)\to \spec(\Oo)$ which is not a smoothing, then there exists a cycle $Z$ with $\chi(\sO_Z)\le 1$, which lifts to $\XX$.
\end{prop}

\subsection{Equivariant deformations in which $Z$ lift}

Let $Z$ be the fundamental cycle for the minimal resolution 
$X\to \spec(R)$. 
We are interested in the deformations of $X$ to which the divisor $Z$ lifts.  Let $R_Z^G$ be the functor defined by 
$$ 
R_Z^G (A) =\Big\{\text{$G$-equivariant~}\XX\to \spec (A)|  \text{the Cartier divisor $\ZZ$ lifting $Z$}
\Big\}/\cong
$$
In the minimal resolution, $H^0(N_Z)=0$, therefore such $\ZZ$ are uniquely determined by $\XX$. This defines a subspace 
$R_Z^G\subset D_X^G$, the $G$-equivariant deformation space of $X$. 

From \cite[Proposition 2.2]{Wahl_Duke}, this functor $R_Z^G$ has a formal versal deformation space.  We have the following result for the space $R_Z^G$.
For the divisor $i: Z\subset X$, we have the exact triangle 
for the tangent complex (which is the dual of the cotangent complex)
\begin{equation}\label{eqn_tangent_Z-1}
   \ttt_{Z}\to \ttt_X\to \ttt_{Z/X}[1],
\end{equation}
where $\ttt_{Z/X}[1]=N_Z$. Therefore taking cohomology we have a morphism 
\begin{equation}\label{eqn_tangent_Z}
T_Z\to i^*T_X=T_X|_{Z}\to N_Z.
\end{equation}

\begin{prop}\label{prop_alpha_G}
    Let $\alpha: H^1(T_Z)^G\to H^1(N_Z)^G$ be the morphism obtained by taking the $G$-equivariant cohomology for the morphism  (\ref{eqn_tangent_Z}). Then 
    \begin{enumerate}
        \item The tangent space of $R_Z^G$ is $\ker(\alpha)$.
        \item The obstruction space of $R_Z^G$ is $\coker(\alpha)=H^1(T_Z^1)^G$.
    \end{enumerate}
\end{prop}
\begin{proof}
    The first statement is from the standard equivariant deformation theory. 

    The second statement is also from the standard deformation and obstruction theory.  Suppose that $A\to \overline{A}$ is a small extension with kernel $(\epsilon)$.  Let $\overline{Z}\subset \overline{X}$ be an element  of $R_Z^G(\overline{A})$.  Lifting the minimal resolution $\overline{X}$ to $X$ over $\spec(A)$, and there exists an obstruction $\ob\in H^1(N_Z)^G$ to lifting $\overline{Z}$ to $X$.  If the image of $\ob$ in $\coker(\alpha)$ is zero, then $\ob=\alpha(\tau)$ for a $\tau\in H^1(T_{\overline{X}})^G$.

    From (\ref{eqn_tangent_Z-1}), there is an exact sequence
    $$
    0\to T_Z\longrightarrow T_X\longrightarrow N_Z\longrightarrow T_Z^1\to 0
    $$
    and 
    $$
    H^1(T_X)\stackrel{\alpha}{\longrightarrow} H^1(N_Z)\longrightarrow H^1(T_Z^1)\to 0.
    $$
    Thus, the cokernel $\coker(\alpha)=H^1(T_Z^1)^G$.
\end{proof}

\subsection{Equivariant deformation of $X$}

Still for the minimal resolution $X\to \spec(R)$, let $B^G\subset D_X^G$ be the subfunctor of the $G$-equivariant deformations of $X$ which blow down to give equivariant deformations of $R$. Similar argument as in \cite{Wahl_Ann} implies that $B^G$ has a formal versal deformation space 
$\overline{B}^G$. Also there is a map 
$$\Phi: B^G\to D_R^G$$
which is finite from \cite[Theorem 3.2]{Wahl_Duke}.

For the normal surface singularity $(\spec(R),0)$, we can write $R=P/I$, where $P$ is a regular local ring of dimension $e$, and $I$ is the ideal. Then the finite group $G$ acts on the ring $P$ and the ideal $I$. 
We denote $T_P$ the tangent space of $\spec(P)$, and after choosing the generators of $\Mm_P$ there is a map 
\begin{equation}\label{eqn_map_TXe}
T_X\to \sO_X^e.
\end{equation}

\begin{prop}\label{prop_beta_map}
    Let $\beta: H^1(T_X)^G\to H^1(\sO_X^e)^G$ be the map obtained by taking $G$-equivariant cohomology on (\ref{eqn_map_TXe}). Then 
    \begin{enumerate}
        \item  The tangent space of $B^G$ is $\ker(\beta)$.
        \item   The obstruction space of $B^G$ is $\coker(\beta)$.
        \item  $\ker(B^G(\mathbf{k}[\epsilon])\to D_R^G(\mathbf{k}[\epsilon]))=H^1_{E}(T_X)^G$.
    \end{enumerate}
\end{prop}
\begin{proof}
    A $G$-equivariant deformation $X^\prime$ of $X\to \spec(R)$ blows down  to give a $G$-equivariant deformation $R^\prime$ of $\spec(R^\prime)\to \spec(A)$ if and only if 
    $H^0(\sO_{X^\prime})\to H^0(\sO_X)$ is $G$-equivariant and surjective.
    Let $\Mm_R=(g_1,\cdots, g_e)$ be the generators. The above is equivalent to the lifting for each $g_i$. The obstruction of the above process is in $H^1(\sO_X)^G$.  Thus, we can perform the same proof as in Theorem \ref{prop_alpha_G} to prove (1), (2).
    (3) is a consequence in \cite[1.8]{Wahl_Ann}.
\end{proof}

\subsection{The proof of Theorem \ref{thm_Wahl-2}}

Recall that the minimally elliptic singularities were introduced by Laufer \cite{Laufer}. A normal singularity $(\spec(R),0)$ is called {\em minimally elliptic} if it is Gorenstein and $h^1(\sO_X)=1$ for the minimal resolution $X\to \spec(R)$. 

Let $Z$ be the fundamental cycle, then $K_X=\sO_X(-Z)$. When the negative self-intersection number $d:=-Z\cdot Z\le 4$, the singularity is lci, and when $-Z\cdot Z\le 3$, it is a hypersurface singularity. When  $d\ge 3$, $R$ is isolated and when $d\ge 2$, $d$ equals to the embedding dimension. 
Minimally elliptic singularities include types {\em simple elliptic singularity; node; cusp; tacnode; triangular} corresponding to the exceptional curves of the minimal resolution. 
We first have a result \cite[Proposition 5.3]{Wahl_Duke}.

\begin{prop}\label{prop_Wahl-Duke5.3}
    Let $X\to (\spec(R),0)$ be the minimal resolution of the  minimally elliptic singularity $(\spec(R),0)$, and $Z$ be the fundamental cycle. Then $H^1(T_Z^1)=0$ if and only if $d=1$ or $2$.  In this case, $Z$ is smoothable. 
\end{prop}

Suppose that we have a finite group $G$ action on the normal surface singularity $(\spec(R),0)$.  The quotient $\spec(R^G,0)$ is still a normal surface singularity.  There is a $G$-equivariant minimal resolution 
$X\to \spec(R)$ such that the quotient 
$\overline{X}\to \spec(R^G)$ is also a partial resolution of $\spec(R^G,0)$.

Recall that  $B^G\subset D_X^G$ is the subspace consisting of the  $G$-equivariant deformations  $\XX\to \spec(A)$ that blow down to give $G$-equivariant deformations of $\spec(R)$. 
Let $R_Z^G\subset D_X^G$ be the subspace consisting of $G$-equivariant deformations such that the exceptional cycle $Z$ lifts $G$-equivariantly. 

\begin{thm}\label{thm_1.2_not_intro}
    Let $X\to (\spec(R),0)$ be the minimal resolution of the  minimally elliptic singularity $(\spec(R),0)$, and $Z$ be the fundamental cycle.
    Suppose that there is a finite group $G$ action on the minimal resolution so that it is $G$-equivariant. Then we have 
    \begin{enumerate}
        \item  $B^G=R_Z^G$;
        \item  We have the map $B^G\to NF^G$, and 
        $\overline{B}^G/W\cong \overline{NF}^G$ when $d\ge 3$. 
        \item  When $d\ge 5$, $NF^G\subset D_R^G$.
    \end{enumerate}
\end{thm}
\begin{proof}
    We generalize the proof in \cite[Theorem 5.6]{Wahl_Duke} to the $G$-equivariant setting. First we have the following equivariant diagram 
    \[
    \xymatrix{
    \ttt_X\ar[d]\ar[r]& \sO_X^{e}\ar[d]\\
    \sO_Z(Z)\ar[r]& \sO_Z^e
    }
    \]
where the left vertical morphism is from (\ref{eqn_tangent_Z}).  The bottom morphism was the one in  \cite[Theorem 5.6]{Wahl_Duke}, but works $G$-equivariantly. We review it here.  Let $\Mm_R$ be the maximal ideal of $R$ with generators 
$(g_1, \cdots, g_e)$, and $\overline{g}_i\in H^0(\sO_Z(Z))^G$. Then the morphism is given by 
$$
\Mm_R\sO_X\to  H^0(\sO_X(-Z))^G\to  H^0(\sO_Z(-Z))^G.
$$
We prove $R_Z^G=B^G$. It suffices to prove that the morphism 
$$
H^1( \sO_Z(Z))^G\to H^1(\sO_Z)^e
$$
is injective, because if it is true, $R_Z^G\to B^G$ would be sujective on the tangent space and injective on the obstruction spaces. Thus, it implies that they must be isomorphic. 

By Serre duality we prove the sujectivity of the morphism 
$$
H^0( \omega_Z(-Z))^e\to H^0(\omega_Z\otimes \sO_Z(-Z))
$$
under the $G$-action. In the minimally elliptic case, we have 
$$
\omega_Z=\omega_X\otimes \sO_Z(Z)\cong \sO_Z
$$
since $\omega_X-\sO_X(-Z)$.  Thus it is enough to prove that 
$$
(H^0( \sO_Z^e)^G\to H^0(\sO_Z(-Z))^G
$$
is surjective. This morphism if given by 
$(\overline{g}_1, \cdots, \overline{g}_e)$, and 
$
(H^0( \sO_Z^e)^G\to H^0(\sO_X(-Z))^G=(\Mm_R)^m$
is surjective for some $m\in \zz_{>0}$. The morphism 
$H^0(\sO_X(-Z))^G\to H^0(\sO_Z(-Z))^G$ is surjective since 
$H^1(\sO_X(-2Z))=0$ from \cite[Theorem 5.6]{Wahl_Duke}.

The proof of  (2), (3) are similar to \cite[Theorem 5.6]{Wahl_Duke} which works $G$-equivariantly. 
\end{proof}

Recall that the types of minimally elliptic singularities are Types 
El (simple elliptic singularity), No (Node), Cu (cusp singularity), 
Ta (Tacnode), Tr (Triangle singularity).  Simple elliptic singularities and cusp singularities are log canonical surface singularities which lie in the semi-log-canonical surfaces in the KSBA moduli spaces. 
We have the following  result.

\begin{cor}
    Suppose that the singularity $(\spec(R),0)$ is of type $\El, \No, \Cu, 
\Ta$, or $\Tr$. Then the fundamental cycle $Z=E=\sum_i F_i$, and 
\begin{enumerate}
    \item $B^G$ and $NF^G$ are smooth of dimension 
    $h^1(ES)^G+h^0(T_E^1)^G$, where 
    $h^1(ES)$ is the number of $G$-equivariant equisingular deformation. 
    \item $\spec(R)$ is the normally flat specialization of a cone.
    \item In the simple elliptic singularity $\El$ case, if the negative self-intersection number $-Z\cdot Z\ge 10$, 
    $NF^G\subset D_R^G$ is an irreducible component. 
\end{enumerate}
\end{cor}
\begin{proof}
    For (1), we need to calculate 
    $\dim(H^1(T_X))^G-\dim(H^1(N_Z))^G$. 
    In this case the divisor $E$ has only simple normal crossing singularities and the equisingular deformations of $R$ is given by 
    $h^1(ES)^G$.
    We have the exact sequence
    $$
    0\to ES\longrightarrow T_X\longrightarrow N_E\longrightarrow T_E^1\to 0.
    $$
    We have $H^0(N_E)=H^1(T_E^1)=0$,  hence 
    $$
    \dim(H^1(T_X))^G-\dim(H^1(N_Z))^G=\dim(H^1(ES))^G+\dim(H^0(T^1_E))^G.
    $$

    (2) is similar to \cite[Corollary 5.7]{Wahl_Duke}.
    For (3), Pinkham \cite{Pinkham} proved that the simple elliptic singularity $(\spec(R),0)$ admits a smoothing if and only if the degree or the negative self-intersection number $1\le -Z\cdot Z\le 9$. \cite{Jiang_2023} studied when such a simple elliptic singularity admits an lci cover by an lci simple elliptic singularity. Thus, when $-Z\cdot Z\ge 10$, the deformations form irreducible components. 
\end{proof}

\section{Equivariant deformation of minimally elliptic singularities}\label{sec_deformation_elliptic_singularity}

\subsection{$G$-equivariant $\omega^*$-constant deformation}\label{subsec_omega_constant}

\subsubsection{}

Let $R$ be a normal, $n$-dimensional Cohen-Macaulary  domain, finite type over $\mathbf{k}$.  Let $\omega:=\omega_R$ be the dualizing module, which is rank one and depth $n$.  Let $G$ be a finite group acting on $\spec(R)$ so that $R^G$ is still Cohen-Macaulay.  We can write $R=P/I$, where $P$ is a regular local ring of dimension $r$, such that there exists a resolution by projective $P$-modules:
$$
0\to F_{r-n}\longrightarrow \cdots \longrightarrow F_1\longrightarrow
P\longrightarrow P/I\to 0.
$$
The group $G$ acts on the resolution, and we have 
$\omega=\coker(F_{r-n-1}^*\to F_{r-n}^*)\cong \Ext_P^{r-n}(R,P)$ which admits a $G$-action. 
The dualizing sheaf can be defined as follows. Let $U\subset \spec(R)$ be the smooth locus, and let $i: U\to \spec(R)$ be the inclusion. Then 
$\omega=i^*\Omega_U^n$.

If the dualizing sheaf $\omega$ is invertible, then the singularity $(\spec(R),0)$ is Gorenstein. Under the $G$-action, the dualizing sheaf $\omega^G:=\omega_{\spec(R^G)}$ of the quotient $\spec(R^G)$ must be $\qq$-Gorenstein. 

\begin{example}\label{example_cusp_simple_Z2-1}
    Let $(\spec(R),0)$ be a type $\El$ singularity with degree $d$. From \cite{Jiang_2023},  
    in the case that $d=8,9$, there exists a $\zz_2$ or $\zz_3$-action such that 
    $\omega^{\zz_2}$ is still Gorenstein and the quotient is a degree $4$ and $3$ simple elliptic singularity. 
\end{example}

\begin{example}\label{example_cusp_simple_Z2-2}
    Here is a cusp singularity studied in \cite[\S 6, Example 1]{Jiang_2022}. 
    Let $(\spec(R),0)$ be a type $\Cu$ singularity whose resolution cycle is given by
\begin{equation}\label{eqn_cusp_graph}
\xymatrix@R=9pt@C=24pt@M=0pt@W=0pt@H=0pt{
  &\overtag{\bullet}{-e_2}{8pt}\dashto[r]
  &&&\overtag{\bullet}{-e_{k-1}}{8pt}\dashto[l]\\
  \lefttag{\bullet}{2-2e_1}{8pt}\lineto[ur]\lineto[dr] &&&&&
  \righttag{\bullet}{2-2e_k}{8pt}\lineto[ul]\lineto[dl]\\
  &\overtag{\bullet}{-e_2}{8pt}\dashto[r]
  &&&\overtag{\bullet}{-e_{k-1}}{8pt}\dashto[l]}
  \end{equation}
where $k\ge 2$, $e_i\ge 2$ and some $e_j>2$. 
The $\zz_2$-quotient-cusp singularity  $(\spec(R^{\zz_2}),0)$ has resolution graph
\begin{equation}\label{eqn_cusp_graph2}
\xymatrix@R=7pt@C=24pt@M=0pt@W=0pt@H=0pt{
  \overtag{\bullet}{-2}{8pt}\lineto[dr] && &&&
  \overtag{\bullet}{-2}{8pt}\lineto[dl]\\
  &\overtag{\bullet}{-e_1}{8pt}\lineto[r]
  &\overtag{\bullet}{-e_2}{8pt}\dashto[r]&\dashto[r]&
  \overtag{\bullet}{-e_k}{8pt}&&\\
  \overtag{\bullet}{-2}{8pt}\lineto[ur] && &&&
  \overtag{\bullet}{-2}{8pt}\lineto[ul]}
  \end{equation}
  The normal singularity $(\spec(R),0)$ is Gorenstein, but $(\spec(R^{\zz_2}),0)$ is $\qq$-Gorenstein with index $2$. 

  There are cusp singularities $(\spec(R),0)$  which admit finite group action such that the quotients $(\spec(R^{G}),0)$  are still cusp singularities which is Gorenstein, see examples in \cite{Jiang_cusp}.
\end{example}

\subsubsection{}

Let $p: \rM\to M$ be a flat surjective morphism of finite type whose fibers are Coehn-Macaulary, of dimension $\ge 2$, with isolated singularities.  Suppose that $G$ acts on $\rM$ so that $G$ preserves the fibers. Then the relative dualizing sheaf $\omega_{\rM/M}$ admits a $G$-action and commutes with base change. 

Consider $\omega^*_{\rM/M}=\Hom(\omega_{\rm/M},\sO)$, the dual. Then    $\omega^*_{\rM/M}$ may not be flat over $M$, and may not commute with base change. 

\subsubsection{}

Let $S$ be a $G$-equivariant Cohen-Macaulay scheme, of dimension $\ge 2$, with only isolated singularities. We have the deformation functor $D_S^G$. For any $G$-equivariant deformation $S^\prime\to \spec(A)$, we have the $\sO_{S^\prime}$-module $\omega_{S^\prime/A}$ flat over $A$ which induces $\omega_S$. 

\begin{defn}\label{defn_omega*_constant}
    A $G$-equivariant deformation  $S^\prime\to \spec(A)$ is 
    $\omega^*$-constant if 
    \begin{enumerate}
        \item $\omega_{S^\prime/A}^*$ reduces $\omega_S^*$ on $S$.
        \item  $\omega_{S^\prime/A}^*$ is $A$-flat. 
    \end{enumerate}
\end{defn}

Let $L_S^G\subset D_S^G$ be the subfunctor of $D_S^G$ of the above $G$-equivariant deformations. 

\begin{prop}(\cite[Proposition 1.7]{Wahl_Math-Ann})
    The subfunctor $L_S^G\subset D_S^G$ has a good deformation theory, and is represented by a closed scheme of the formal versal deformation space $D_S^G$. 
\end{prop}
\begin{proof}
    As in the proof \cite[Proposition 1.7]{Wahl_Math-Ann}, from Schlessinger's criterion \cite{Schlessinger}, Wahl proved that a fibered sum of deformations with flat $\omega^*$ still has flat $\omega^*$. This is because relative dualizing sheaf $\omega$ of a fibered sum over $A\otimes A^\prime$ is isomorphic to the fibered sum of relative $\omega_A$ and $\omega_{A^\prime}$. Then taking dual yields the $\omega^*$.  This process works $G$-equivariantly. 
\end{proof}

The following result gives the criterion of $\omega^*$-constant deformation. 
\begin{prop}\label{prop_Wahl-1.9_Math-Ann}(\cite[Proposition 1.9]{Wahl_Math-Ann})
    Let $R$ be a Cohen-Macaulay domain of dimension $\ge 2$ with isolated singularity.  Let $\spec(\Rr)\to \spec(A)$ is a $G$-equivariant deformation of $R$ with $A$ a discrete valuation ring. Then $\omega^*_{\Rr/A}\cong \omega^*_{\Rr}$, and the $G$-equivariant deformation is $\omega^*$-constant if and only if the depth $\omega^*_{\Rr}\ge 3$. 
\end{prop}

\subsection{$G$-equivariant $\omega^*$-constant deformations}\label{subsec_omega_constant}

Let us study some $G$-equivariant $\omega^*$-constant deformations. 
Let $R=P/I$ be a flat $A$-algebra, with fibers Cohen-Macaulay domains of dimension $n$, where $P$ is a smooth $A$-algebra of relative dimension $r$. From above we have 
\begin{equation}\label{eqn_omega*}
    \omega^*_{R/A}=\ker(F_{r-n}\otimes R\to F_{r-n-1}\otimes R).
\end{equation}

Let us recall the $2$-dimensional cyclic quotient singularities. 
Let $R=\mathbf{k}[x,y]^H$, where $H\subset \GL(2)$ is a cyclic group of order $n$, and the action is given by 
$\mat{cc}\eta & 0\\ 0& \eta^q\rix$, where $\eta$ is a primitive $n$-th root of $1$, and $0< q<n$, $(q,n)=1$.  This is the $n/q$ cyclic quotient singularity.  One writes $n/q$ as a continued fraction expansion
$$
\frac{n}{q}=b_1-\frac{1}{b_2-\frac{1}{b_3-\cdots-\frac{1}{b_k}}}. 
$$
Then the resolution cycle of the cyclic quotient singularity is given by 

\begin{equation}\label{eqn_cusp_graph3}
\xymatrix@R=7pt@C=24pt@M=0pt@W=0pt@H=0pt{
  &\overtag{\bullet}{-b_1}{8pt}\lineto[r]
  &\overtag{\bullet}{-b_2}{8pt}\dashto[r]&\dashto[r]&
  \overtag{\bullet}{-b_k}{8pt}&&}
  \end{equation}

Let $\sum_j (b_j-2)=m-2=e-3$,  then $m$ is the multiplicity, and $e$ is the embedding dimension.  In the case $m=4$, Riemenschneider proved that the versal deformation space of this $n/q$ cyclic quotient singularity has tangent space dimension 
$$
t=\sum_j (b_j-1)+1.
$$
The deformation space contains two smoothing components, of dimension 
$t-1$ (the Artin component), and $t-3$ (the non-Artin component $NA$). These components can be explicitly written down, see \cite[(2.5.1), (2.5.2)]{Wahl_Math-Ann}.

Let $(\spec(R),0)$ be a type $n/q$ cyclic quotient singularity, then our finite group $G$ acts on the singularity $(\spec(R),0)$ with the quotient singularity $(\spec(R^G,0))$.
\cite[Theorem 2.6]{Wahl_Math-Ann} proved that if $\spec(R)$ is a cyclic quotient singularity of multiplicity $4$, then $L_R^G$ is smooth, and is represented by the non-Artin component $NA$.

Here is an example of $\omega^*$-constant smoothing of Wahl singularity. 
Let $R=\mathbf{k}[x,y]^H$, where $H$ is a cyclic group of order 
$n$ such that $(\spec(R),0)$ is a type $n/(n-1)$ singularity (which is called $A_{n-1}$-singularity).  It is an lci singularity. 
Let $G=\mu_n$ be the cyclic group of order $n$ again with primitive root of unity $\eta$.  Let $\eta$ act on $P=\mathbf{k}[x,y,z]$ by 
$$
x\mapsto \eta x;  \quad y\mapsto \eta^q y; \quad z\mapsto \eta^{n-1}z.
$$
Then $t=xz-y^n\in P^G$. $P^G$ is Cohen-Macaulay, and the flat morphism $\mathbf{k}[t]\to P^G$ gives a deformation 
$$
\spec(P^G)\to \spec(\mathbf{k}[t])=\aaa_{\mathbf{k}}^1.
$$
We have 
$R\cong \mathbf{k}[x,y,z]/(xz-y^n)$, and it admits a smoothing 
$$\spec(\Rr)=\spec(\mathbf{k}[x,y,z,t]/(t-xz-y^n))\to \spec(\mathbf{k}[t]).$$ 
The cyclic group $G$ acts on $\spec(\Rr)$ as above.  
We calculate that $R^G$ is a cyclic quotient singularity 
of type $n^2/(nq-1)$ for $(q,n)=1$. Then we get

\begin{prop}\label{prop_Wahl_singulairty_G}
    The $G$-equivariant deformation 
    $\spec(\Rr)\to \spec(\mathbf{k}[t])$ induces the deformation
    $\spec(P^G)\to \spec(\mathbf{k}[t])$.
    The smoothing $\spec(\Rr)\to \spec(\mathbf{k}[t])$ is $\omega^*$-constant, and not on the Artin component. 
\end{prop}
\begin{proof}
    The first statement is from the construction.  The second is from \cite[Theorem 2.7]{Wahl_Math-Ann}. 
\end{proof}

\begin{rmk}
    The above $\omega^*$-constant smoothing can be generalized to type $T$-singularities $\spec(\mathbf{k}[x,y])/\mu_{n^2\cdot s}$, where $\eta$ is a $n^2\cdot s$-th root of unity acting by 
    $\mat{cc}\eta&0\\0&\eta^q\rix$ for $0< q<n$ and $(q,n)=1$.

    The $G$-equivariant deformation 
    $\spec(\Rr)\to \spec(\mathbf{k}[t])$  is the $\qq$-Gorenstein deformation  in the KSBA moduli space in \cite{Kollar-Shepherd-Barron}.
\end{rmk}

\subsection{Partial resolution, blowing down and minimally elliptic singularities}\label{subsec_partial_resolution_G}

Let $X\to \spec(R)$ be the $G$-equivariant minimal resolution of a minimally elliptic surface singularity $(\spec(R),0)$. 
The dualizing sheaf $\omega_X\cong \sO_X(-Z)$, where $Z$ is the fundamental cycle and the degree $d=-Z\cdot Z$ is the multiplicity of $\spec(R)$.

The main
Theorem \ref{thm_Wahl-2} implies that a $G$-equivariant deformation of $X$ blows down to a $G$-equivariant deformation of $\spec(R)$ if and only the fundamental cycle $Z$ lifts to the deformation.  If $Z$ is smoothable, such deformations of $X$ yield an equidimensional subvariety $M_m$ of the deformation space $D_R^G$ which consists of all the $G$-equivariant elliptic deformations of $\spec(R)$ with the same $-Z\cdot Z$. The generic point in $M_m$ corresponds to simple elliptic singularities $\El(d)$ such that $-Z\cdot Z=d$.

We fix $F_1,\cdots, F_s$ to be the disjoint exceptional rational configurations in $X$, and let 
$$f: X\to Y$$
be the $G$-equivariant blowing-down.  The blowing down gives rational singularities $Y_i=\spec(R_i)$ on $Y$ at the singularities $p_1, \cdots, p_s$.  The induced map $Y\to \spec(R)$ is called the ``partial resolution" of $\spec(R)$. 
Then we have 
$I=f_*\sO_X(-Z)\cong \omega_Y$ since $f_*\omega_X=\omega_Y$. 

The following is from \cite[Lemma 3.5]{Wahl_Math-Ann}.
\begin{enumerate}
    \item $f_*\omega_X^*=\omega_Y^*$, and $\omega_Y^*=I^*=\Hom(I,\sO_Y)$.
    \item  Apply $\Hom(I,-)$ to the exact sequence 
    $$
    0\to I\longrightarrow \sO_Y\longrightarrow \sO_Y/I\to 0
    $$
    and we get 
    the exact sequence
    $$
    0\to \sO_Y\longrightarrow I^*\longrightarrow \Hom(I, \sO_Y/I)\to 0
    $$
    since $\Hom(I,I)=\sO_Y, \Ext^1(I,I)=0$. 
    \item $f_*\sO_Z(Z)\cong \Hom(I,\sO_Y/I)=\Hom(I/I^2,\sO_Y/I)$, which is obtained using (2) and the exact sequence 
    $$
    0\to \sO_X\longrightarrow \omega_X^*\longrightarrow \sO_Z(Z)\to 0.
    $$
\end{enumerate}

Now we are interested in the $G$-equivariant deformations $Y^\prime\to \spec(A)$ to which $I$ lifts. This means there is an $I^\prime\subset \sO_{Y^\prime}$ which is flat over $A$ and induces $I$ on $\sO_Y$.  All of these are $G$-equivariant. 
The deformation space of lifting $I$ is 
$$
H^0(\Hom(I/I^2, \sO_Y/I))=H^0(Y, f_*\sO_Z(Z))=H^0(X, \sO_Z(Z))=0
$$
from \cite[Lemma 3.5]{Wahl_Math-Ann}.  Thus, the lifting is unique if there is one.  The obstruction space is $\Ext^1(I/I^2, \sO_Y/I)$ and from Leray spectral sequence it breaks up to 
\begin{multline*}
  0\to H^1(\Hom(I/I^2, \sO_Y/I))\longrightarrow 
\Ext^1(I/I^2, \sO_Y/I) \\
\longrightarrow H^0(\sE xt^1(I/I^2, \sO_Y/I))\to 0, 
\end{multline*}
where the first space is to patch globally, and the last space gives the local  obstruction to lift $I$. 

Thus, we have the following result from \cite[Theorem 3.7]{Wahl_Math-Ann}.

\begin{thm}\label{thm_Wahl_3.7}(\cite[Theorem 3.7]{Wahl_Math-Ann})
    Let $Y^\prime\to \spec(A)$  be a $G$-equivariant deformation of $Y$. Then 
    \begin{enumerate}
        \item  If the ideal $I^\prime\subset \sO_{Y^\prime}$ lifts $I$, then $Y^\prime$ blows down to a $G$-equivariant deformation of $\spec(R)$. Also we have 
        $I^\prime\cong \omega_{Y^\prime/A}$.
        \item  On the other hand, if $[Y^\prime]\in D_Y^G(A)$ is a $G$-equivariant deformation, and $Y^\prime$ blows down, then 
        $I^\prime$ lifts to $I$.
    \end{enumerate}
\end{thm}

\subsubsection{$\omega_Y^*$-constant deformation}

Let $D_Y^G$ be the $G$-equivariant deformation space of $Y$, and 
$$B_Y^G\subset D_Y^G$$
be the subfunctor of $G$-equivariant deformations of $Y$ to which $\omega_Y^*$ lifts and blows down to $G$-equivariant deformations of $\spec(R)$.

We have the following main result. 

\begin{thm}\label{thm_G_BYG-DYi}
    Let $X\to Y\to \spec(R)$ be a $G$-equivariant partial resolution 
    of a minimally elliptic singularity, with rational configurations $\{F_i\}$ on $X$ and $G$-equivariantly blown down to rational singularities $Y_i$ on $Y$. If the fundamental cycle $Z$ is reduced off the support of the $F_i$'s, then the natural map
    $B_Y^G\to \prod D_{Y_i}^G$ is smooth. 
\end{thm}

\begin{proof}
    For the partial resolution $Y\to \spec(R)$, the obstruction 
    $$
    H^2(Y,T_Y)\cong H^1(T_Y^1)=0
    $$
    since $T_Y^1$ only supports on isolated rational singularities. Thus, 
    $H^2(Y,T_Y)\cong H^1(T_Y^1)^G=0$.
    This means that 
    $$
    D_Y^G\to \prod D_{Y_i}^G
    $$
    is smooth.  So to prove $B_Y^G\to \prod D_{Y_i}^G$ is smooth, it suffices to show 
    \begin{equation}\label{eqn_surjective1}
        H^1(Y,T_Y)\to H^1(\Hom(I/I^2, \sO_Y/I))
    \end{equation}
    is surjective under the $G$-action.  Since $f_*T_X=T_Y$, it suffices to show 
    \begin{equation}\label{eqn_surjective2}
        H^1(f_*T_X)\to H^1(f_*\sO_Z(Z))
    \end{equation}
    is surjective under the $G$-action. 
    Look at the inclusion $i: \hookrightarrow X$, and the exact sequence
    $$
    I_Z/I_Z^2\to \rightarrow i^*\Omega_X\rightarrow \Omega_Z\to 0,
    $$
    apply the functor $R\Hom(-,\sO_Z)$ to the above exact sequence we get 
    $$
    T_X\rightarrow \sO_Z(Z)\rightarrow \Ext^1(\Omega_Z, \sO_Z)\to 0.
    $$
    The space $\Ext^1(\Omega_Z, \sO_Z)$ is supported at the singular points of $Z$, and $f_*T_X)\to f_*\sO_Z(Z)$ has cokernel with finite support. Thus, (\ref{eqn_surjective2}) is surjective. 
\end{proof}

For the partial resolution $X\to Y\to \spec(R)$, let $E_1, \cdots, E_s$ be the exceptional curves on $X$ which intersect none of $\{F_i\}$. Then $f: X\to Y$ is isomorphic on $E_j$, which we still denote them as $E_j$ in $Y$. Let 
$$
J:=f_*\sO_X(-Z+\sum_j E_j)=I(\sum_j E_j),
$$
then $J$ is a sheaf of ideals on $Y$.  Let 
$$
LB_Y^G\subset B^G_Y
$$
be the $G$-equivariant deformations to which $E_j$ lifts, and equivalently, $LB_Y^G$ are those $G$-equivariant deformations to which the ideals $\sO(E_j)$ and $J$ lift. 

\begin{thm}\label{thm_LBG}
    Suppose that $Z$ is reduced off the support of $\{F_i\}$. Then 
    the map $LB_Y^G\to \prod D_{Y_i}^G$ is smooth. 
\end{thm}
\begin{proof}
    Let $S_Y$ be the sheaf of dual of the differentials which preserve $E_j$, then there is an exact sequence
    $$
    0\to S_Y\rightarrow T_Y\rightarrow \oplus_j E_j.
    $$
    We have the map $H^1(T_Y)\to \oplus_j H^1(N_{E_j})$ is surjective and $G$-equivariant. This implies that 
    $$LD_Y^G\to \prod D_{Y_i}^G$$
    is smooth, where $LD_Y^G$ is the $G$-equivariant deformation functor to which $E_j$ lift.  We actually want to obtain 
    $LB_Y^G$, and this is the $G$-equivariant deformation such that 
    $J\subset \sO_Y$ also lifts. The obstruction theory requires that 
    \begin{equation}\label{eqn_obstruction_surjective}
        H^1(S_Y)\to H^1(\Hom(J/J^2, \sO_Y/J))
    \end{equation}
     is surjective under the $G$-action. To prove (\ref{eqn_obstruction_surjective}), note that $f_*S_X=S_Y$, where $S_X$ is defined by 
     $$
     0\to S_X\rightarrow T_X\rightarrow \oplus_j N_{E_j}\to 0.
     $$
     Let $Z_1:=Z\setminus \{E_j\}$, then (\ref{eqn_obstruction_surjective}) becomes 
     $H^1(f_*S_X)\to H^1(f_*\sO_{Z_1}(Z_1))$.
     The sheaf $S_X$ equals to $T_X$ on all but a finite number of points of $Z_1$, and one concludes that 
     the cokernel $f_*S_X\to f_*\sO_{Z_1}(Z_1)$ has finite support.  Thus, (\ref{eqn_obstruction_surjective}) is surjective under the $G$-action. 
\end{proof}

\subsubsection{Equivariant blowing down $B^G_Y\to D_R^G$}\label{subsec_G_blow_down}

Recall the $G$-equivariant deformation space $B_Y^G$, and we like to study the equivariant blow-down map
$$
\Psi: B_Y^G\to D_R^G.
$$

The following result \cite[Proposition 4.2]{Wahl_Math-Ann} is automatically true when we consider the $G$-equivariant deformations. 

\begin{prop}\label{prop_4.2_Wahl-Math-Ann}(\cite[Proposition 4.2]{Wahl_Math-Ann})
    Let $X\to Y\to \spec(R)$ be the $G$-equivariant minimal resolution and partial resolution of a minimally elliptic singularity.  Suppose that $Y$ only has rational singularities, and $Y^\prime$ is the generic fiber of a small $G$-equivariant deformation of $Y$ such that it $G$-equivariantly  blows down. Let $X^\prime\to Y^\prime$ be the $G$-equivariant minimal resolution of singularities, then 
    $X^\prime\to \Gamma(\sO_{X^\prime})=\spec(R^\prime)$ is the $G$-equivariant minimal resolution. 
\end{prop}

Let $\{Y_i\}$ be the rational singularities of $Y$, and 
$\{Y_j^\prime\}$ be the rational singularities of $Y^\prime$. Then we have the following formula of the self-intersection of the fundamental cycles.
\begin{equation}\label{eqn_Z-Z_number}
-Z\cdot Z-\sum \ell(Y_i)=-Z^\prime\cdot Z^\prime-\sum \ell(Y_j^\prime),
\end{equation}
where $\ell(Y_i)$ is defined as follows. 
The rational singularity $Y_i=\spec(R_i)$, and 
$\ell(Y_i)=\Ext^1(\omega_{R_i}, R_i)$.
If $I_i\subset R_i$ is the canonical ideal $(I_i\cong \omega_{R_i})$, we have 
$$
\ell(Y_i)=\dim(\Ext^1(I_i, R_i/I_i))=\dim(\Ext^2(R_i/I_i, R_i/I_i)).
$$
Also if $R_i$ is a cyclic quotient singularity, we have 
$$\ell(Y_i)=\max(0, \mult R_i-3).$$
Suppose that $R_i$ admits a finite group $G$-action, then we have 
$$\ell(Y^G_i)=\max(0, \mult R^G_i-3).$$

\begin{thm}\label{thm_minimal_resolution_G_Z}
    Let $X\to Y\to \spec(R)$ be the $G$-equivariant minimal resolution and partial resolution of a minimally elliptic singularity.  Suppose that 
    \begin{enumerate}
        \item $\{Y_i\}$ are cyclic quotient  singularities, either of multiplicity $4$ or as in Proposition \ref{prop_Wahl_singulairty_G}.
        \item The fundamental cycle $Z$ is reduced off the support of $\{F_i\}$.
    \end{enumerate} 
    Let $m_i=\mult(Y_i^G)$. Then $\Psi(B_Y^G)\subset D_R^G$ is irreducible, and the generic point corresponds to a simple elliptic singularity of degree 
    $-Z\cdot Z-\sum (m_i-3)$.
\end{thm}
\begin{proof}
    Theorem \ref{thm_G_BYG-DYi} implies that $B_Y^G\to \prod D_{Y_i}^G$ is smooth. Thus, from Proposition \ref{prop_Wahl_singulairty_G}, there exists a one-parameter $G$-equivariant deformation in $B_Y^G$
    whose generic fiber $Y^\prime$ is nonsingular.  Then $Y^\prime$ is the $G$-minimal resolution of a minimally elliptic singularity and from (\ref{eqn_Z-Z_number}) we have 
    $$
    -Z^\prime\cdot Z^\prime=-Z\cdot Z-\sum \ell(Y_i^\prime)=
    -Z\cdot Z-\sum (m_i-3).
    $$
    Since the  fundamental cycle $Z$ is reduced off the support of $\{F_i\}$, $Z^\prime$ is reduced. 
    Thus, $Y^\prime$ can be deformed so that $Z^\prime$ can $G$-equivariantly blow down to a cone over an elliptic curve with the same $-Z^\prime\cdot Z^\prime$.
\end{proof}

It is well-known from \cite{Laufer} that a minimally elliptic singularity of degree $d\le 4$ is lci. When $d=5$, it is not lci, but is given by  the Pfaffians. The deformation of all of these cases  does not have higher obstructions, see \cite{Jiang_2021}.
We let 
$$
\scM\to M
$$
be the versal universal family of $G$-equivariant deformations of $\spec(R)$.  Suppose that $\spec(R)$ is a minimally elliptic singularity. 
Define
$$
\Delta_i=\text{Closure of~} \{s\in M | \scM_s \text{~has an~} \El(i) \text{~singularity~}\}.
$$
Then 
$$
\bigcup_{i}\Delta_i=\{s\in M | \scM_s \text{~has an $G$-equivariant~} \El\text{~singularity~}\}.
$$

\begin{prop}\label{prop_G_deformation_property}
    Suppose that $\scM\to M$ is the $G$-equivariant versal deformation space of a minimally elliptic singularity $\spec(R)$ of multiplicity $m\le 4$.  Then 
    \begin{enumerate}
        \item $M$ is smooth.
        \item $\Delta_m$ is of pure dimension $\dim(M)+m-10$.
        \item $\Psi(B_Y^G)$ is an irreducible component of $\Delta_l$, where $l=-Z\cdot Z-\sum (m_i-3)$.
    \end{enumerate}
\end{prop}
\begin{proof}
    In the case that $m\le 4$,  the minimally elliptic singularity is lci.  If the minimally elliptic singularity $\spec(R)$ admits a $G$-action, then it definitely admits a $G$-equivariant deformation. (1) is true, i.e.,  $M$ is smooth, and the dimension calculations (2), (3) are from \cite[Theorem 4.8]{Wahl_Math-Ann}. 
\end{proof}

\subsection{Lci minimally elliptic singularities}\label{subsec_cusp_simple_G_allpication}

Let $\spec(R)$ be a minimally elliptic singularity with multiplicity $m$. When $m\le 4$, the singularity is lci.  
We are mainly interested in $\El$, and $\Cu$ singularities. 
The hypersurface $\El$ singularities were given in \cite{Looijenga2} corresponding to $E_6, E_7, E_8$ types.
The local equation of lci  cusp singularities were given in \cite{Nakamura}.
A $\Cu$ singularity $(\spec(R),0)$ is a local complete intersection ($\lci$) cusp if and only if $(\spec(R),0)$ is one of the following:
    \begin{enumerate}
        \item $T_{p,q,r}: x^p+y^q+z^r-xyz=0$ with $\frac{1}{p}+\frac{1}{q}+\frac{1}{r}<1$,
        \item $\prod_{p,q,r,s}: x^p+w^r=yz,\quad  y^q+z^s=xw$ with $(\frac{1}{p}+\frac{1}{r})(\frac{1}{q}+\frac{1}{s})<1$,
    \end{enumerate}
    where $p,q,r$ and $p,q,r,s$ are integers greater than $1$ and the cusp point  is chosen to be the origin. 
If a finite group $G$ acts on the regular ring $P$ so that $R=P/I$ where $I$ is the ideal generated by the local equations, then these singularities admit $G$-equivariant deformations. 

Other type of minimally elliptic singularities contain Triangle singularities, $\Ta$, and $\No$ singularities, see \cite[\S 5.1]{Wahl_Math-Ann} for their resolution graphs. 

Suppose that the $G$-action on $\spec(R)$ induces a Gorenstein quotient singularity
$(\spec(R^G),0)$ which is still a minimally elliptic singularity. For example,  a $\El$ singularity of degree  $4$ or  $3$ admits a $\zz_2$, $\zz_3$ action whose quotient is a degree $8$ or $9$ $\El$-singularity, see \cite{Jiang_2023}.  $\Cu$-singularities admit finite group action whose quotients are still cusp singularities, see \cite{Pinkham}, \cite{Jiang_cusp}
for more examples. 

First we have
\begin{lem}\label{lem_minimal_rational_G}
    Let $(\spec(R),0)$ be a Gorenstein surface singularity, together with a finite group action. Then 
    \begin{enumerate}
        \item suppose that $(\spec(R),0)$ is a rational singularity, then the quotient $(\spec(R^G),0)$ is still rational, and a $G$-equivariant $\omega^*$-constant deformation of $(\spec(R),0)$ induces $\omega^*$-constant deformation of $(\spec(R^G),0)$.
        \item Suppose that $(\spec(R),0)$ and the quotient $(\spec(R^G),0)$ are both minimally elliptic singularities, then any $G$-equivariant deformation of $(\spec(R),0)$ induces a deformation of $(\spec(R^G),0)$. 
    \end{enumerate}
\end{lem}
\begin{proof}
    Let $\spec(\Rr)\to \spec(A)$  be a $G$-equivariant deformation with $A$ a discrete valuation ring. Then $\Rr^G$ is also a Cohen-Macaulay ring, $\qq$-Gorenstein since $G$ is finite.  Thus, $\spec(\Rr^G)\to \spec(A)$ gives a deformation of $(\spec(R^G),0)$.

    For (1), let $f: X\to \spec(R)$ be the $G$-equivariant minimal resolution.
    Then $f^G: X/G\to \spec(R^G)$ is a partial resolution. Let 
    $$\tilde{f}: \tilde{X}\to X/G\to \spec(R^G)$$
    be the minimal resolution.  We write $\tilde{f}=f^G\circ \pi$ where 
    $\pi: X\to X/G$. Then 
    $$R^1\tilde{f}_{*}\sO_{\tilde{X}}
    =R^1(f^G\circ \pi)_{*}\sO_{\tilde{X}}=R^1f^G_*(\sO_{X/G})=0$$
    since $R^1f_*(\sO_X)=0$. Thus, $\spec(R^G)$ is a rational singularity. 

    Suppose that $\spec(\Rr)\to \spec(A)$ is $\omega^*$-constant. From 
    Proposition \ref{prop_Wahl-1.9_Math-Ann}, $\omega^*_{\Rr/A}\cong \omega^*_{\Rr}$, and the $G$-equivariant deformation is $\omega^*$-constant if and only if the depth $\omega^*_{\Rr}\ge 3$. 
    Let $t\in A$ be a regular parameter. Then the depth $(\omega_{\Rr}^*)\ge 3$ if and only if 
    the depth $(\omega_{\Rr}^*/t \omega_{\Rr}^*)\ge 2$, if and only if 
    $\omega_{\Rr}^*/t \omega_{\Rr}^*=\omega_{R}^*$.  This is also true for the quotient since 
    $\omega_{\Rr^G}^*/t \omega_{\Rr^G}^*=\omega_{R^G}^*$ since 
    $\spec(\Rr^G)\to \spec(A)$ is a deformation of $(\spec(R^G),0)$.  
    Thus, it is $\omega^*$-constant. 

    (2) is from the deformation theory of minimally elliptic singularities.
    Let $\ll_R$ be the cotangent complex of $R$.  Then the $G$-equivariant deformation space of $\spec(R)$ is given by 
    $\Ext^1(\ll_R, R)^G$ which is not zero since the quotient $\spec(R^G)$ is still a minimally elliptic singularity. 
\end{proof}    

Generalizing \cite[Theorem 5.4, Theorem 5.6]{Wahl_Math-Ann} we have the following $G$-equivariant deformation components of minimally elliptic singularities.  For the readers' convenience, we list the dual graph of known minimally elliptic singularities. 
\begin{equation}
    \begin{array}{ll}
      \El(d):   & \text{simple elliptic singularity} \\
       \Cu(d_1, \cdots, d_r):  &  \text{cusp with resolution cycle of $\pp^1$'s}  \\
       D(b_1, b_2, b_3): & \text{Dolgachev singularity},
    \end{array}
\end{equation}
where the dual graph of Dolgachev singularity is given by\\
\begin{equation}\label{eqn_Dolgachev}
\xymatrix@R=7pt@C=24pt@M=0pt@W=0pt@H=0pt{
  &&\overtag{\bullet}{-b_2}{8pt}\lineto[dd] && &&&\\ \\
  &\undertag{\bullet}{-b_1}{1pt}\lineto[r]
  &\undertag{\bullet}{-1}{1pt}\lineto[r]&
  \undertag{\bullet}{-b_3}{1pt}&&}
  \end{equation} \\ \\
and the singularities have three cases according to the minimal resolution:
$\No(d)=D(2,3,6+d)$ has Node: 
\begin{figure}[htbp]
  \centering
  \includegraphics[width=20mm]{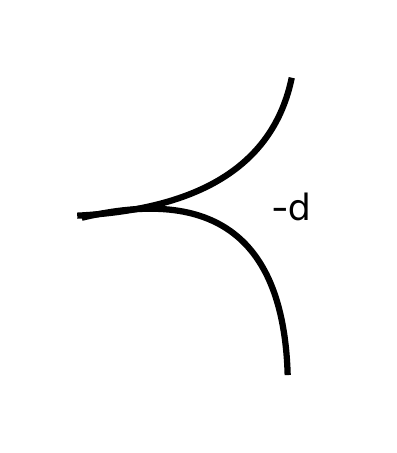}      
  \caption{Node}
  \label{fig-No-config}
\end{figure}

\begin{figure}[htbp]
  \centering
  \includegraphics[width=20mm]{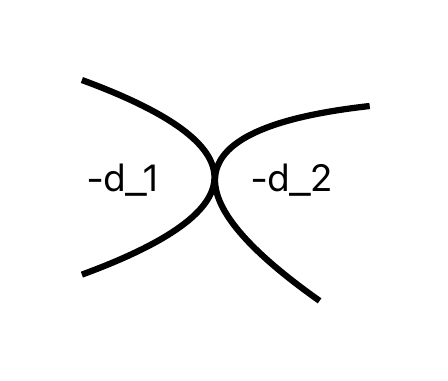}      
  \caption{Tacnode}
  \label{fig-Ta-config}
\end{figure}

$\Tr(d_1, d_2)=D(2,d_1+2,d_2+2)$ has the following singularity, see Figure \ref{fig-Ta-config}.

\begin{figure}[htbp]
  \centering
  \includegraphics[width=20mm]{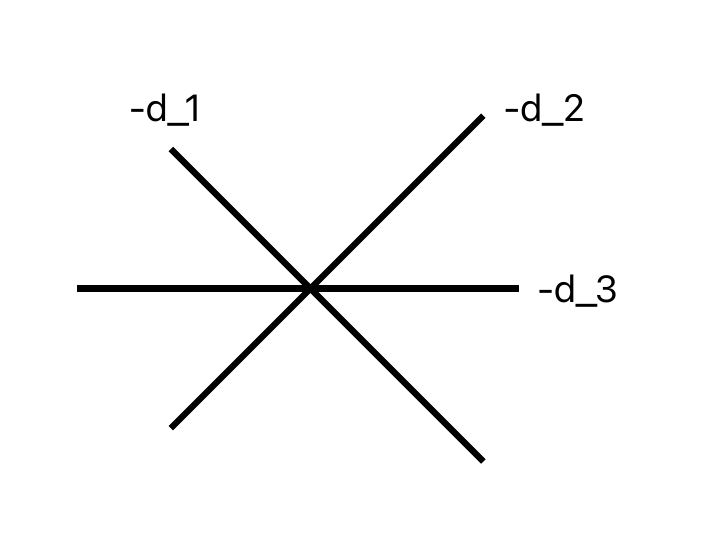}      
  \caption{Triangle}
  \label{fig-Tr-config}
\end{figure}

$\Tr(d_1,d_2,d_3)=D(d_1+1,d_2+1,d_3+1)$ hashas the following singularity, see Figure \ref{fig-Tr-config}. 

We let singularities $S_1\longrightarrow S_2$ represent $S_1$ $G$-equivariantly deforms to $S_2$, if there is a $G$-equivariant one-parameter deformation family with special fiber $S_1$, and generic fiber $S_2$.
This also means $S_2$ is $G$-equivariantly adjacent to $S_1$. 
We define the following.
\begin{equation}\label{eqn_deformations_G_Cu_El}
    \begin{array}{ll}
      (1)   &  \Cu(d_1, d_2,\cdots, d_r)\longrightarrow \Cu(d_1+d_2-2,d_3,\cdots, d_r) (r>2) \\
         & \Cu(d_1, d_2)\longrightarrow \Cu(d_1+d_2-4)\\
         & \Cu(d)\longrightarrow \El(d) \\
      (2)& \No(d) \longrightarrow \Cu(d) \\
      (3)& \Ta(d_1,d_2)\longrightarrow \Cu(d_1, d_2)\\
         & \Ta(d_1,d_2)\longrightarrow \No(d_1+d_2-4)\\
      (4)& \Tr(d_1,d_2,d_3)\longrightarrow \Cu(d_1, d_2, d_3)\\
         & \Tr(d_1,d_2,d_3)\longrightarrow \Ta(d_1+d_2-2,d_3)
    \end{array}
\end{equation}
In all of the above we assume that there is a finite group $G$-action on the deformations so that they are $G$-equivariant deformations. 

\begin{thm}\label{thm_G_equisingular_minimal}
    Let $(\spec(R),0)$ be a minimally elliptic singularity and $X\to \spec(R)$ the $G$-equivariant minimal resolution. Let $Z$ be the fundamental cycle and let $d:=-Z\cdot Z$ which is given by the following formula for the type of minimally elliptic singularities.
    \begin{equation}\label{eqn_ZZ_G_Cu_El}
    \begin{array}{ll}
      \Cu(d_1, d_2,\cdots, d_r):   & -Z\cdot Z=\begin{cases}
          \sum (d_i-2) & r>1\\
          d & r=1.
      \end{cases}  \\
      \No(d): & -Z\cdot Z=d \\
    \Ta(d_1,d_2): & -Z\cdot Z=d_1+d_2-4 \\
    \Tr(d_1,d_2,d_3): & -Z\cdot Z=d_1+d_2+d_3-6
    \end{array}
\end{equation}
    Then the $G$-equivariant deformations (\ref{eqn_deformations_G_Cu_El}) (which are $-Z\cdot Z$-constant deformations) form an irreducible component of the deformation space $D_R^G$. 
\end{thm}
\begin{proof}
    The results are from the general result Theorem \ref{thm_Wahl-2} and Theorem \ref{thm_1.2_not_intro} by considering  $Z$ lifting the $G$-equivariant deformations. 
\end{proof}

\begin{example}
    Let $(\spec(R_1),0)$ be a $\Cu(2,5)$-singularity. From \cite[Lemma 2.5]{Nakamura}, it is a hypersurface singularity given by 
    $$x^3+y^3+z^5-xyz=0.$$
    Let $(\spec(R_2),0)$ be a $\Cu(3)$-singularity. It is a hypersurface singularity given by 
    $$x^3+y^3+z^4-xyz=0.$$
    Let $\mu_3$ act on $R_1, R_2$ by:
    $$
    x\mapsto \eta x, \quad  y\mapsto \eta^2 y, \quad z\mapsto z.
    $$
    From Theorem \ref{thm_G_equisingular_minimal}, $R_1$ $\mu_3$-equivariantly deforms to $R_2$. This can be seen in their infinitesimal deformation space $T^1_{R_i}:=T^1_{\spec(R_i)}$ for $i=1, 2$. 
    Both $\Cu(2,5)$ and $\Cu(3)$  have $-C^2=3$, where $C=C_1+\cdots +C_r$ is the cycle of $\pp^1$'s. 
    From \cite[Theorem 5.3]{Nakamura2}, using the notations there, we have 
    $(b_0, b_1, b_2)=(3,3,5)$ for $R_1$; and $(3,3,4)$ for $R_2$.  So
    $\dim T^1_{R_1}=3+3+5+1=12$, and $\dim T^1_{R_2}=3+3+4+1=11$. Therefore, the deformation space of $R_2$ is contained in the deformation space of $R_1$. 

    Let us look at the quotient. From \cite[Example 1, \S 2.3]{Jiang_cusp}, the quotient of $R_1$ is still a $\Cu(3, 2,2,2,2,2)$-singularity which is a hypersurface $(\spec(R_1^{\mu_3}),0)$
    $$x^2+y^3+z^{12}-xyz=0.$$
    For the cusp $\Cu(3)$, the cyclic group $\mu_3$ is the torsion subgroup of the first integral homology of its link, then from \cite[Theorem 1]{Pinkham}, its quotient is the cusp $\Cu(2,2,3)$, which is also a hypersurface $(\spec(R_2^{\mu_3}),0)$
     $$x^2+y^3+z^{9}-xyz=0.$$
     From \cite[Theorem 5.3]{Nakamura2}, $\dim(T^1_{R_2^G})=2+3+9+1=15$, and
    from \cite[Theorem 5.3]{Nakamura2}, $\dim(T^1_{R_1^G})=2+3+12+1=18$. 
   Thus, $R_1^G$ deforms to $R_2^G$. 
\end{example}

The following result deals with $G$-equivariant deformations that involve lowering the degree of $-Z\cdot Z$.  We are mainly interested in the $\El$ and $\Cu$ singularities.  For other type of minimally elliptic singularities, see \cite[Theorem 5.6]{Wahl_Math-Ann}.

\begin{thm}\label{thm_G_singular_minimal_deformation}
    Let $(\spec(R),0)$ be a minimally elliptic singularity and $X\to \spec(R)$ be the minimal resolution of $G$-equivariant. Let $Z$ be the fundamental cycle and let $d:=-Z\cdot Z$ which is given by (\ref{eqn_ZZ_G_Cu_El}).

    Let $\Upsilon$ be a configuration of a rational singularity given by 
\begin{equation}\label{eqn_configuration_rational}
\xymatrix@R=7pt@C=24pt@M=0pt@W=0pt@H=0pt{
  &\undertag{\bullet}{-b_1}{1pt}\lineto[r]
  &\undertag{\bullet}{-b_2}{1pt}\dashto[rr]&&
  \undertag{\bullet}{-b_t}{1pt}&&}
  \end{equation} \\ 
in Proposition \ref{prop_Wahl_singulairty_G}.
Then we have the following $G$-equivariant deformations 
    \begin{equation}\label{eqn_deformations_G_Cu_El-2}
    \begin{array}{ll}
      (1)   &  \Cu(d_1,b_1,\cdots, b_t, d_2,\cdots, d_r)\longrightarrow \Cu(d_1+d_2-1,d_3,\cdots, d_r) (r>2) \\
       (2) & \Cu(d_1,b_1,\cdots, b_t, d_2)\longrightarrow \Cu(d_1+d_2-3)\\
       (3)  & \Cu(d, b_1,\cdots, b_t)\longrightarrow \El(d-1).
    \end{array}
\end{equation}
\end{thm}
\begin{proof}
Let 
$f: X\to Y$ be the $G$-equivariant blowing down so that $Y\to \spec(R)$ is a partial resolution. From Theorem \ref{thm_G_BYG-DYi}, the $G$-equivariant deformations $Y^\prime$ of $Y$ lie in 
$\Psi(B_Y^G)$.
Let $K\subset Y$ be the subscheme defined by $f_*\sO_X(-Z)$.  This $K$ is reduced. Using $LB_Y^G\subset B^G_Y$ and Theorem \ref{thm_LBG}, there is a 
$G$-equivariant smoothing $Y^\prime$ of $Y$ and a reduced subscheme $K^\prime$ consisting of $E_3,\cdots, E_r$ and a lifting $C^\prime$ of $f(E_1\cup E_2)$.  
A further deformation of $Y^\prime$, if necessary, can leave $E_3, \cdots, E_r$ unchanged and smooth $C^\prime$ to give the desired configuration in the theorem. 
\end{proof}

\begin{rmk}
    In Theorem \ref{thm_G_singular_minimal_deformation}, the deformation quotient of the simple elliptic and cusp singularities will induce the deformation of their quotients.
\end{rmk}



\subsection*{}


\begin{thebibliography}{12}  
 

\bibitem{Alper_moduli} Jarod Alper, \newblock  Stacks and Moduli, https://sites.math.washington.edu/~jarod/moduli.pdf.


 



\bibitem{Jacob_equi} J.J. Byszewski, Cohomological aspects of equivariant deformation theory, PhD thesis.

\bibitem{Engel} P. Engel,  \newblock A proof of Looijenga's conjecture via integral-affine geometry, 
{\em Journal of Differential Geometry},  109(3): 467-495, 2018. 


\bibitem{GHK15} M. Gross, P. Hacking and S. Keel,  \newblock  Mirror symmetry for log Calabi-Yau surfaces I, {\em Publ. Math. Inst. Hautes Etudes Sci.}, 122:65-168, (2015). 



\bibitem{Jiang_2021}Y. Jiang, \newblock A note on higher obstruction spaces for surface singularities, preprint, arXiv:2112.10679. 

\bibitem{Jiang_2022}Y. Jiang, \newblock The virtual fundamental class for the moduli space of surfaces of general type, arXiv:2206.00575.
\bibitem{Jiang_2023}Y. Jiang, \newblock  Smoothing of surface singularities via equivariant smoothing of lci covers,  preprint, arXiv:2309.16562.  
\bibitem{Jiang_cusp}Y. Jiang, \newblock  Equivariant smoothing of cusp singularities,  preprint, arXiv:2302.00637.

\bibitem{Kollar-Shepherd-Barron} J. Koll\'ar and N. I. Shepherd-Barron, \newblock Threefolds and deformations of surface singularities, 
{\em Invent. Math.}, 91, 299-338 (1998).

\bibitem{Kresch_Tsch} Andrew Kresch, Yuri Tschinkel, Birational geometry of Deligne-Mumford stacks, 	arXiv:2312.14061.

\bibitem{Laufer} H. Laufer, \newblock On minimally elliptic singularities, {\em Amer. J. Math.}, Vol. 99, No. 6 (1977), 1257-1295.

\bibitem{Looijenga}E. Looijenga, \newblock  Rational surfaces with an anticanonical cycle, {\em Ann. of Math.} (2), 114 (2):267-322, 1981.

\bibitem{Looijenga2}E. Looijenga, \newblock, On the semi-universal deformation of a simple-elliptic hypersurface singularity: Part I: Unimodularity, {\em Topology}
Vol. 16, Iss. 3, 1977,  257-262.

\bibitem{Nakamura} I. Nakamura, \newblock  Inoue-Hirzebruch surfaces and a duality of hyperbolic unimodular singularities. I. {\em Math. Ann.} 252, 221-235 (1980). 
\bibitem{Nakamura2} I. Nakamura, \newblock Infinitesimal Deformations of Cusp Singularities, {\em Hokkaido Mathematical Journal},  Vol. 15(1984), p. 21-46.


\bibitem{Pinkham}H. C. Pinkham, \newblock  Automorphisms of cusps and Inoue-Hirzebruch surfaces, {\em Compositio Math.} 52, No.3 (1984), 299-313. 


\bibitem{Schlessinger}M. Schlessinger, \newblock Functors of Artin rings, {\em Trans. Amer. Math. Soc.} 130, 208-222 (1968).


\bibitem{Wahl_Ann} J. Wahl, \newblock Equisingular deformations of normal surface singularities, I, {\em Ann. Math.}, Vol. 104, No.2, 325-356 (1976). 
\bibitem{Wahl_Duke} J. Wahl, \newblock Simultaneous resolution and discriminantal loci, {\em  Duke Math. J.},  Vol. 46, No.2,  341-275 (1979).
\bibitem{Wahl_Math-Ann} J. Wahl, \newblock Elliptic deformation of minimally elliptic singularities, {\em Math. Ann.}, 253, 241-262 (1980).
\end{thebibliography}
\end{document}